
\font\vlbf=cmbx10 scaled 1440
\font \lbf=cmbx10 scaled 1200
\baselineskip 7mm plus 0pt minus 0pt
\hfill{}

\vskip 3cm
\hfill{q-alg/9508007}
\centerline{\vlbf On the Classification of quantum qroup stuctures on the
group GL(2)}

\vskip 2cm

\centerline{\bf A. Aghamohammadi ${}^{1,2,*}$, M. Khorrami ${}^{1,3,4}$,
A. Shariati ${}^{1,4}$}

\vskip 1cm
{\it
\noindent $^1$ Institute for Studies in Theoretical Physics and
	    Mathematics, P. O. Box  5746, Tehran 19395, Iran.

\vskip -2mm\noindent $^2$ Department of Physics, Alzahra University,
	      Tehran 19834, Iran.

\vskip -2mm\noindent $^3$ Department of Physics, Tehran University,
	     North-Kargar Ave. Tehran, Iran.

\vskip -2mm\noindent $^4$ Institute for Advanced Studies in Basic Sciences,
	     P. O. Box 159, Gava Zang, Zanjan 45195, Iran.

\vskip -2mm\noindent $^*$ E-Mail: mohamadi@irearn.bitnet}

\vskip 2cm
\centerline{\lbf Abstract}

All quantum group structures on the group GL(2) are classified.
It is shown that there are only two such structures, the well known
quatum groups GL$_{qp}$(2) and GL$_{hh'}$(2).

\noindent
\vfil\break
Although the group GL(2) is not a complicated one, and its dimension is not
so high, it seems that there are still some ambiguities in the
classification of the quantum group structures on GL(2).

There are two well-known deformations of this group [1]: the two parametric
Drinfeld-Jimbo deformation GL$_{q,p}$(2) [2], and the two parametric
Jordanian deformation GL$_{h,h'}$(2) [3]. It has been shown [4] that these
two parametric deformations are related to each other by a singular
transformation (contraction).
It has been proved that, up to isomorphisms, there are only two quantum
group structures on GL(2) with a central quantum determinant,
GL$_q$(2) and GL$_h$(2) [5].
Recently Kupershmidt has tried to classify all
quantum group structures on the group GL(2) [6]. He has concluded that there
are three different quantum group structures on GL(2): GL$_{q,p}$(2),
GL$_{h,h'}$(2), and the new structure GL$_{h,q}$(2). Here we will show that
this new structure is, in fact, not new: by a suitable regular
transformation, it falls in one of the two well-known structures
GL$_{q,p}$(2) or GL$_{h,h'}$(2), depending on the deformation parameters.
So, in fact there are only two non-equivalent quantum group structures on
GL(2).

To begin, we introduce the quantum-plane relations of GL$_{q,p}$(2) and
GL$_{h,h'}$(2). For GL$_{q,p}$(2) we have
$$\eqalign{xy&=qyx,\cr \xi^2&=0,\cr \eta^2&=0,\cr \eta\xi&=-p\xi\eta .\cr}
  \eqno{(1)}$$
And for GL$_{h,h'}$(2)
$$\eqalign{[x,y]&=hy^2,\cr \xi^2&=h'\xi\eta ,\cr \eta^2&=0,\cr \eta\xi &=
  -\xi\eta .\cr}\eqno{(2)}$$
Here $x$ and $y$ are the coordinates of the quantum plane, $\xi =dx$, and
$\eta =dy$. The relations of the new quantum plane introduced in ref. [6]
are
$$\eqalign{xy&=yx,\cr \xi^2&=h'_0\xi\eta ,\cr \eta^2&=r_0\xi\eta ,\cr
  \eta\xi &=-p_0\xi\eta .\cr}\eqno{(3)}$$
In fact, this is not the parametrization used in ref. [6]. We have used this
parametrization for later convenience.

As it has been noted in ref. [6], one can use a transformation to make
$r=0$. (Applying the transformation changes $(h'_0,r_0,p_0)$ to $(h',r,p)$.)
The point is that, by a similar transformation, one can make both $h'$ and
$r$ equal to zero. So that this new structure is, in fact, the old
GL$_{q,p}$(2) with $q=0$. There is, however, one exeption: if $(p_0-1)^2=
4h_0'r_0$, one cannot make both $h'_0$ and $r_0$ equal to zero. Instead, one
can make $r$ equal to zero, and $p$ equal to one. But this is GL$_{h,h'}$(2)
with $h=0$. Here we present the above-mentioned transformations.

\noindent{\bf I} Case $(p_0-1)^2\ne 4h_0'r_0$.

\noindent Using the linear transformation
$$\pmatrix{\bar\xi\cr\bar\eta\cr}=\pmatrix{1&t_1\cr 1&t_2\cr}
   \pmatrix{\xi\cr\eta\cr},\eqno{(4)}$$
where $t_i$'s are the roots of the equation
$$r_0t^2+(1-p_0)t+h'_0=0,\eqno{(5)}$$
one can see that
$$\bar\xi^2=\bar\eta^2=0,\eqno{(6)}$$
and
$$\eqalign{\bar\xi\bar\eta&=(h'_0+t_2-p_0t_1+r_0t_1t_2)\xi\eta ,\cr
 \bar\eta\bar\xi&=(h'_0+t_1-p_0t_2+r_0t_1t_2)\xi\eta .\cr}\eqno{(7)}$$
Note that, $(p_0-1)^2\ne 4h_0'r_0$ guarantees that $t_1\ne t_2$, and this
shows that the transformation (4) is nonsingular. Also note that the
coefficients of $\xi\eta$ in (7) cannot be both zero; otherwise, from (6)
and (7), one would conclude that $\xi^2=\eta^2=\xi\eta =\eta\xi =0$. So, in
this case, the structure is equivalent to GL$_{q=1,p}$(2), where
$$p=-{{h'_0+t_1-p_0t_2+r_0t_1t_2}\over{h'_0+t_2-p_0t_1+r_0t_1t_2}}.
 \eqno{(8)}$$

\noindent{\bf II} Case $(p_0-1)^2=4h_0'r_0$.

\noindent In this case, the transformation (4) is singular. So we use
$$\pmatrix{\bar\xi\cr\bar\eta\cr}=\pmatrix{1&0\cr 1&t\cr}
   \pmatrix{\xi\cr\eta\cr},\eqno{(9)}$$
where $t$ satisfies (5), that is
$$t={{p_0-1}\over{2r_0}}={{2h'_0}\over{p_0-1}}.\eqno{(10)}$$
Now we have
$$\eqalign{\bar\eta^2&=0,\cr\bar\xi^2&={{h'_0}\over{h'_0+t}}
  \bar\xi\bar\eta ,\cr\bar\eta\bar\xi&=-\bar\xi\bar\eta .\cr}\eqno{(11)}$$
We see that, in this case, the sructure is equivalent to GL$_{h=0,h'}$(2)
with
$$h'={{h'_0}\over{h'_0+t}}.\eqno{(12)}$$

These linear transformations on the quantum plane induce similarity
transformations on the quantum matrix. So we have concluded that there are
exactly two distinct deformations of GL(2).
\vfil\break
\noindent{\lbf References}
\item{[1]} H. Ewen, O. Ogievetsky, \& J. Wess; Lett. Math. Phys. {\bf 22}
(1991) 297.
\item{[2]} Yu. I. Manin; {\it Topics in noncommutative geometry}; Princeton,
NJ: Princeton University Press.
\item{[3]} A. Aghamohammadi; Mod. Phys. Lett. {\bf A8} (1993) 2607.
\item{[4]} A. Aghamohammadi, M. Khorrami, \& A. Shariati; J. Phys. {\bf A28}
(1995) L225.
\item{[5]} B. A. Kupershmidt; J. Phys. {\bf A25} (1992) L1239.
\item{[6]} B. A. Kupershmidt; J. Phys. {\bf A27} (1994) L47.
\end